\documentclass[12pt]{article}
\usepackage{amsmath}
\usepackage{amssymb}
\usepackage{tabularx}
\usepackage{enumerate}
\usepackage{mathrsfs}
\usepackage{amssymb, upref, color}
\usepackage[dvips]{graphicx}
\newtheorem{thm}{Theorem}[section]
\newtheorem{lem}[thm]{Lemma}

\newtheorem{rmk}{Remark}[section]
\newtheorem{defi}{Definition}[section]
\newtheorem{pppp}{Proof}

\newcommand{\qed}{\hspace{1em}\mbox{\raisebox{0.65ex}{\fbox{}}}}

\numberwithin{equation}{section}

\newcommand{\be}{\begin{equation}}
\newcommand{\ee}{\end{equation}}
\newcommand\bes{\begin{eqnarray}} \newcommand\ees{\end{eqnarray}}
\newcommand{\bess}{\begin{eqnarray*}}
\newcommand{\eess}{\end{eqnarray*}}

\newcommand{\R}{\mathbb{R}}
\newcommand{\bpf}{{\bf Proof:\ \ }}
\newcommand{\epf}{\mbox{}\hfill $\Box$}

\begin{document}
\thispagestyle{empty}

\title{Spatial spreading model and dynamics of West Nile virus in birds and mosquitoes with free boundary\thanks{The work is partially supported by the
NSFC of China (Grant No. 11371311 and 11171267), the High-End Talent Plan of Yangzhou University and NSERC and CIHR of Canada.}}
\date{\empty}

\author{Zhigui Lin$^a$ and Huaiping Zhu$^{b}$\\
{\small $^a$School of Mathematical Science, Yangzhou University, Yangzhou 225002, China}\\
{\small $^b$Laboratory of Mathematical Parallel Systems (LAMPS)}\\
 {\small Department of Mathematics and Statistics}\\
 {\small York University, Toronto, ON, M3J 1P3, Canada}
}

 \maketitle

\begin{quote}
\noindent {\bf Abstract.}
In this paper, a reaction-diffusion system is proposed to model the spatial spreading of West Nile virus in vector mosquitoes and host birds in North America. Transmission  dynamics are based on a simplified model involving mosquitoes and birds, and the free boundary is introduced to model and explore the expanding front of the infected region. The spatial-temporal risk index $R_0^F(t)$, which involves regional characteristic and time, is defined for the simplified reaction-diffusion model with the free boundary to compare with other related threshold values, including the usual basic reproduction number $R_0$.  Sufficient conditions for the virus to vanish or to spread are given.  Our results suggest that the virus will be in a scenario of vanishing if $R_0\leq 1$, and will spread to the whole region
if $R_{0}^F(t_0)\geq 1$ for some $t_0\geq 0$, while if $R^F_0(0)<1<R_0$, the spreading or vanishing of the virus depends on the initial number of infected individuals, the area of the infected region, the diffusion rate and other factors. Moreover, some remarks on the basic reproduction numbers and the spreading speeds are presented and compared.

\end{quote}

\section{Introduction}

\begin{figure}[htbp] \centering
\includegraphics[height=3.5in]{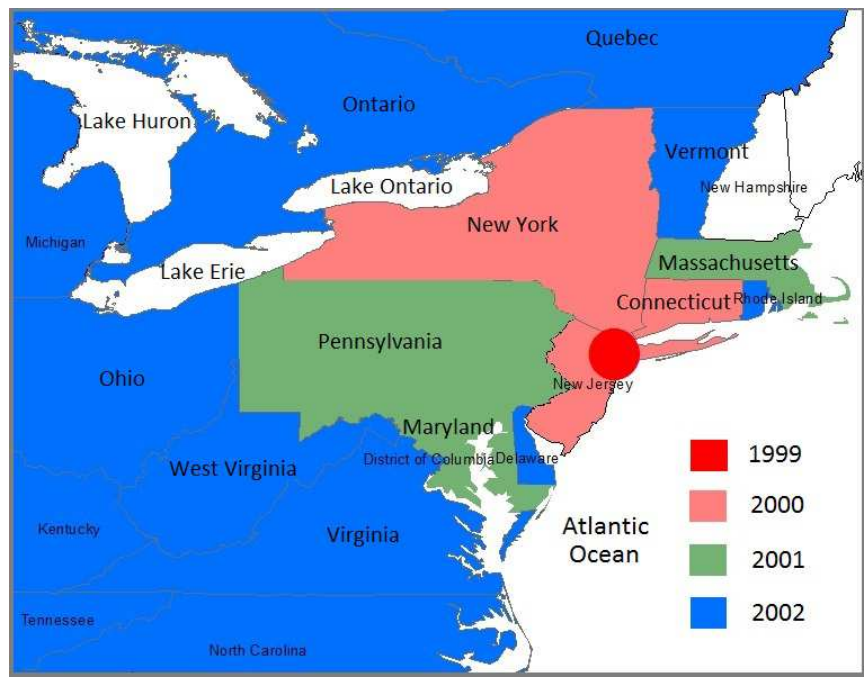}
   \caption{The spatial spreading of West Nile virus from New York city to its neighboring states 
   from 1999 to 2002. }\label{WNV-CDC}
\end{figure}

West Nile virus (WNv) is an arthropod-borne flavivirus that cause the epidemics
of febrile illness and sporadic encephalitis. It is a typical mosquito-borne disease 
with {\it culex} mosquitoes as vectors and birds as hosts of the virus.
The virus arrived and became endemic 
for the first time in North America in the summer of 1999 in New York City, since then it has kept spreading to its neighboring states. In 2012, the States experienced the largest outbreak of the virus and CDC received reports of 5674 cases of human infection  (\cite{CDC}). As of October 2014, a total of 47 states and the District of Columbia have reported WNv activities in USA \cite{CDCWNV}.
For the case in Canada, the virus moved further west and north, and arrived and caused local endemic in southern Ontario in 2001. In 2002, mosquitoes, birds and horses 
in other provinces including Qu\'ebec, Nova Scotia, Manitoba, and Saskatchewan tested positive for the virus. In 2003, WNv activity was also reported in New Brunswick and Alberta (\cite{phac}). As shown in Fig. \ref{WNV-CDC}, since its first arrival in New York City in 1999, the virus has quickly spread across almost the whole continent of North America.

There have been extensive modeling studies for the virus and most of the work focus 
on the temporal transmission dynamics of the virus between vector mosquitoes and host birds. The available compartmental models for WNv focus more on the existence and stability of equilibria (disease free and endemic equilibria), the temporal transmission dynamics are usually characterized and presented in terms of the so called
basic reproduction number. The related results provided theoretical frame work for developing public health strategies for prevention and control of the virus,  see Wonham et al. \cite{Wonham04}, Gustavo et. al. \cite{Gustavo05}, Bowman et. al. \cite{Bowman} and Abdelrazec et. al. \cite{Ahmed14}  and references therein.

Though the compartment models play an important role in understanding the disease dynamics, there has been increasing interest and need in understanding the spatial spreading processes of WNv. The spatial spreading of WNv is much more complex,
it involves the demographics of both vector mosquitoes and host population (birds, horses, humans etc.), and it is closely related to the movement of both vectors and hosts, and the incidence mechanisms over the expanding front 
where mosquitoes bite hosts to pass the virus to cause new infections.
Lewis et al. \cite{Lewis05} initiated the investigation of the temporal-spatial spreading of the diseases considering the movement of birds and mosquitoes.
Their reaction-diffusion model was developed from a temporal model for WNv by Wonham et al. \cite{Wonham04}, with the diffusion terms describing the movement of birds and mosquitoes.  Under moderate assumptions on the cooperative nature of cross-infection dynamics, Lewis et al \cite{Lewis05} proved the existence of traveling waves and calculated the spatial spread rate of infection in a simplified version of the reaction-diffusion model. Results of comparison theorem are used to show that the spread rate of the simplified model may provide an upper bound for the spread rate of a more realistic and complex version of the model.  Liu et. al. \cite{LiuMBE06} studied the directional dispersal of birds and its impact on spreading of the virus.  Maidana and Yang \cite{Maidana09} proposed a spatial model to analyze the WNv propagation across the USA, and studied the traveling wave solutions of the model to determine the speed of disease dissemination. The wave speed was obtained as a function of the model parameters for the purpose of accessing control strategies. The propagation of WNv from New York City to California state was established as a consequence of the diffusion and advection movements of birds. Moreover, their results showed that mosquito movement does not play an important role in the disease dissemination, while bird advection becomes an important factor for lower mosquito biting rates.

Even though the existence of traveling wave gives an estimate of the speed of the spatial spread wave of the virus, it is the asymptotical wave speed that usually gives an approximation of the progressive spreading speed of the virus transmission, and it does not really reflect the spread of the virus in the early stage of the spatial expanding of the infection to larger area.
It is natural to model the spatial spreading of the virus by using a free-boundary, that is, at the boundary front of an infected area, the virus expands and pushes forward to induce further spatial spreading till the whole region or area become endemic.

In the process of temporal-spatial spreading, what makes WNv or vector-borne diseases special is that at the spreading front, infected mosquitoes (vectors) and  birds (hosts), or both can cause new infection and push the (free) boundary of the infected area forward. Therefore we define 
three types of free boundaries to reflect the three cases of  conditions on the interfaces. But for the analysis, in this paper, we will focus only on the case of  WNv positive birds found on the interface. It will be interesting to study and compare with the other two cases of boundary conditions and their roles in determining the basic reproduction numbers and spreading speed.  We leave them as future work.

In this paper, based on the available temporal-spatial modeling studies of WNv, we will establish and study a reaction-diffusion model with free boundary to explore the temporal-spatial transmission of the virus, where the population of the vector mosquitoes is described by a system for the susceptible and infected classes mosquitoes while the dynamics of the host birds is described by an SIS model, the expanding front is expressed by a free boundary which models the spatial expanding of the infection (infected area).  The spatial-temporal risk index $R_0^F(t)$ will be defined for the simplified model with the free boundary to compare with other related threshold values, including the usual basic reproduction number ($R_0$) and other thresholds ($R_0^N$ and $R_0^D$) related to the reaction-diffusion model. Initial values of infected vector mosquitoes and birds, the area of the initial infected region, the diffusion rates and other factors will be combined to develop sufficient conditions for the virus to vanish or to become spatially endemic.

\section{Model Formulation}

For WNv, its transmission involves vector mosquitoes and host birds. As in the temporal models in \cite{Wonham04, Gustavo05} and \cite{Bowman}, we
classify the vector mosquitoes and host birds in the following subgroups:
\begin{itemize}
\item
susceptible mosquitoes with the number $V_s(t)$, infected mosquitoes with the number $V_i(t)$;
 \item
susceptible birds with the number $H_s(t)$ and infected birds with the number $H_i(t)$.
\end{itemize}
As in \cite{Wonham04}, we  assume that the infected birds recover with no immunity to the virus if they survive the infection, therefore we will not have a recovered class, and infected birds become susceptible once they have recovered. If we do not consider the spatial spreading of the virus, adopted from the models in \cite{Wonham04, Gustavo05, Bowman},  an ODE model describing the temporal transmission of the virus reads
\begin{equation}\label{model1}
\left \{
\begin{array}{rcl}
\dfrac{dV_s}{dt}&=&( V_s+(1-q) V_i)G(V_s,V_i)
     -\dfrac{\beta_vV_sH_i}{N_h}-d_v V_s,\\[2ex]
\dfrac{dV_i}{dt}&=& qV_iG(V_s,V_i) + \dfrac{\beta_vV_sH_i}{N_h}
-d_v V_i,\\[2ex]
\dfrac{dH_s}{dt}&=&r_h(H_s+H_i) -\dfrac{\beta_h V_i H_s}{N_h}-d_h H_s+\gamma_hH_i,\\
\dfrac{dH_i}{dt}&=&\dfrac{\beta_h V_i H_s}{N_h}-d_hH_i -\gamma_h H_i,
\end{array} \right.
\end{equation}
 where $G(V_s, V_i)$ is the per capita reproduction rate of the adult vector mosquitoes which can be
 taken as  $r_v\Big(1-\dfrac{V_s+V_i}{K_v}\Big)$, or simply a constant $r_v$ (so called recruitment rate);  $d_v$ is the natural death rate of mosquitoes; $\beta_v$ is the contact transmission rate of hosts
 to vectors; $r_h$ is the recruitment rate of host birds;  $d_h$ is the natural death rate of birds; $\beta_h$ is the  contact transmission rate of the virus from mosquitoes  to birds and $\gamma_h$ is the recovery
 rate of birds recovering from the infection. The parameter  $0<q\ll 1$ measures the vertical transmission rate of the virus in culex mosquitoes \cite{Dohm2002}.

When the 
mosquitoes and birds are in different spatial
locations, the standard method of including the spatial movement 
consists in the introduction of the diffusion terms. We start with one-dimensional case: $-\infty<x<\infty$, thus based on available temporal model for WNv, a spatially extended version of the vector-host model for WNv can be described by
\begin{equation}\label{model2}
\left \{
\begin{array}{rcl}
\dfrac{\partial V_s}{\partial t}-D_v\dfrac{\partial^2V_s}{\partial x^2}&=&(V_s+(1-q)V_i)G(V_s,V_i)
     -\dfrac{\beta_vV_sH_i}{N_h}-d_v V_s,\\[2ex]
\dfrac{\partial V_i}{\partial t}-D_v \dfrac{\partial^2V_i}{\partial x^2}&=& qV_iG(V_s, V_i) + \dfrac{\beta_vV_sH_i}{N_h} -d_v V_i,\\[2ex]
\dfrac{\partial H_s}{\partial t}-D_h \dfrac{\partial^2H_s}{\partial x^2}&=&r_h(H_s+H_i) -\dfrac{\beta_h V_i H_s}{N_h}-d_h H_s+\gamma_hH_i
      ,\\[2ex]
\dfrac{\partial H_i}{\partial t}-D_h \dfrac{\partial^2H_i}{\partial x^2}&=&\dfrac{\beta_h V_i H_s}{N_h}-d_hH_i -\gamma_h H_i
\end{array} \right.
\end{equation}
for $-\infty<x<\infty$ and $t>0$, where the unknowns $V_s(x,t), V_i(x,t), H_s(x,t)$ and $H_i(x,t)$ are the densities of their respective class in the location $x$ at time $t$, $D_v,\, D_h$ represent the diffusion rates for the vector mosquitoes and host birds, respectively, therefore we can assume that $0<D_v\ll D_h$.

For simplicity, we start with considering the case $G(V_s, V_i)=r_v$, and let $r_v=d_v$, $r_h=d_h$.
In other words, we assume that both the density of vector mosquitoes and that of hosts all remain constants. Furthermore, we assume that the density of vector mosquitoes and that of hosts are initially constant in space, then system (\ref{model2})
implies that $(V_s+V_i)(x,t)$ and $(H_s+H_i)(x,t)$ remain constant in space for all time.
Let $N_v^*=V_s+V_i$ and $N_h^*=H_s+H_i$, then the above system can be simplified to \begin{equation}\label{model3}
\left \{
\begin{array}{rcl}
\dfrac{\partial V_i}{\partial t}-D_v \dfrac{\partial^2V_i}{\partial x^2}&=&\dfrac{\beta_v(N_v^*-V_i)H_i}{N_h^*} -r_v(1-q) V_i,\\[2ex]
\dfrac{\partial H_i}{\partial t}-D_h \dfrac{\partial^2H_i}{\partial x^2}&=&\dfrac{\beta_h V_i (N^*_h-H_i)}{N_h^*}-(d_h+\gamma_h) H_i.
\end{array} \right.
\end{equation}

This research is devoted to understanding the spatial transmission mechanisms of WNv, therefore we will pay more attention to the changing of infected domain and consider a vector-host epidemic model with a free boundary, which describes the spreading frontier of the virus in space.

Usually public health units in North America use three criteria
to decide whether the area is infected by the WNv:
\begin{itemize}
\item
Criterion 1. Found WNv positive vector mosquitoes only.

Following the first arrival of the virus in 1999, many public health units in Canada and regions in USA have been running the mosquito surveillance program, which has been successful in alerting the endemic situation of the virus in the heath units.

 \item
Criterion 2. Found WNv positive birds.

Usually when  dead birds were found and tested WNv positive, it would be a firm sign of the activities of the virus in the region.
In southern Ontario, Canada, dead birds were collected for viral test for the period from 2002 till 2006, the test for birds stopped when there were only few reported  human infection cases.

\item
Criterion 3. Both WNv positive mosquitoes and birds are found, with reported WNv human cases.

\end{itemize}
 The above three cases or related criteria have been
 common practice in public health units. In general, lab tests for both birds and mosquitoes are used in regions of Canada to confirm the endemic of the virus.

In this paper, we will consider the second case. Assume that the mosquitoes and birds migrate in the whole region represented by
$(-\infty, \infty)$, and at time $t$,
 WNv positive  birds were found only in the region represented by $g(t)<x<h(t)$,
there is no infected birds in the rest of the region.

As in \cite{LIN}, the length of the expanding distance $h(t+\Delta t)-h(t)$ is assumed
to be proportional to diffusion mediated gradient of $H_i$, leading to
$$h(t + \Delta t)-h(t)\approx -\mu D_h\frac{\partial H_i}{\partial x}\Delta t.$$
Letting $\Delta t\to 0$, we then obtain the condition on the right interface (free boundary)
$$H_i(h(t),t)=0,\quad -\mu D_h\frac{\partial H_i}{\partial
x}(h(t),t)=h'(t).$$
Similarly, the conditions on the left interface
(free boundary) are
$$H_i(g(t),t)=0,\quad -\mu D_h\frac{\partial H_i}{\partial
x}(g(t),t)=g'(t).$$
 In such a case, we have the problem for
$V_i(x,t)$ and $H_i(x,t)$ with free boundaries $x=g(t)$ and $x=h(t)$ as follows,
\begin{equation}\label{a3}
\left \{
\begin{array}{lll}
\dfrac{\partial V_i}{\partial t}-D_v \dfrac{\partial^2V_i}{\partial x^2}=\dfrac{\beta_v(N_v^*-V_i)H_i}{N_h^*}
-r_v(1-q) V_i,\; & g(t)<x<h(t),\, t>0, \\
\dfrac{\partial H_i}{\partial t}-D_h \dfrac{\partial^2H_i}{\partial x^2}=\dfrac{\beta_h V_i (N^*_h-H_i)}{N_h^*}
-(d_h+\gamma_h) H_i,\; & g(t)<x<h(t), \, t>0,\\
V_i(x,t)=H_i(x, t)=0,&x=g(t)\, \textrm{or}\, x=h(t),\, t>0,\\
g(0)=-h_0,\; g'(t)=-\mu D_h \frac{\partial H_i}{\partial x}(g(t), t), & t>0, \\
 h(0)=h_0, \; h'(t)=-\mu D_h \frac{\partial H_i}{\partial x}(h(t), t), & t>0,\\
V_i(x,0)=V_{i0}(x),\; H_i(x,0)=H_{i0}(x), &-h_0\leq x\leq h_0,
\end{array} \right.
\end{equation}
where $x=g(t), x=h(t)$ are the moving left and right
boundaries to be determined,  $h_0,\ d$ and $\mu $ are positive constants, and the initial functions
$V_{i0}$ and $H_{i0}$ are nonnegative and satisfy
\begin{eqnarray}
\left\{
\begin{array}{ll}
V_{i0}\in C^2([-h_0, h_0]),\, V_{i0}(\pm h_0)=0\, \textrm{and} \ 0\leq V_{i0}(x)\leq N_v^*,\, x\in (-h_0, h_0), \\
H_{i0}\in C^2([-h_0, h_0]),\, H_{i0}(\pm h_0)=0\, \textrm{and} \ 0<H_{i0}(x)\leq N_h^*,\, x\in (-h_0, h_0).
\end{array} \right.
\label{Ae}
\end{eqnarray}
Ecologically, the model means that beyond the free boundaries $x=g(t)$ and $x=h(t)$, there are only susceptible host birds, no birds carrying the virus.

For Criterion 1, the conditions on the interfaces are
$$V_i(g(t),t)=0,\quad -\mu D_v\frac{\partial V_i}{\partial
x}(g(t),t)=g'(t),$$
$$V_i(h(t),t)=0,\quad -\mu D_v\frac{\partial V_i}{\partial
x}(h(t),t)=h'(t).$$

For Criterion 3, the conditions on the interfaces are
$$(V_i+H_i)(g(t),t)=0,\quad -\mu (D_v\frac{\partial V_i}{\partial
x}+D_h\frac{\partial H_i}{\partial
x})(g(t),t)=g'(t),$$
$$(V_i+H_i)(h(t),t)=0,\quad -\mu (D_v\frac{\partial V_i}{\partial
x}+D_h\frac{\partial H_i}{\partial
x})(h(t),t)=h'(t).$$
In this case, we assume that $V_i+H_i>0$ in $(g(t), h(t))$, which implies that in the area there are mosquitoes or birds infected with the virus.

In the rest of the paper, we will only consider problem (\ref{a3}) for Criterion 2 and give the properties of the solution and the free boundaries, similar discussions can be done for Criteria 1 and 3.

\section{Existence and uniqueness}

In this section, we first present the following basic results on the existence and uniqueness of the spreading model \eqref{a3} with the initial conditions \eqref{Ae}.
\begin{thm}\label{exist}
The following hold:
\begin{itemize}
 \item[{\rm (i)}] For any given $V_{i0}, H_{i0}$ satisfying $(\ref{Ae})$, problem $(\ref{a3})$
admits a unique global solution $(V_i, H_i, h, g)$ with $V_i, H_i\in C^{(1+\alpha)/2, 1+\alpha}(\overline{D}_{T})$ and
$h, g\in C^{1+\alpha/2}([0,T])$
 for any $\alpha \in (0,1)$ and $T>0$, where  $D_{T}=\{(x,t)\in \mathbf R^2:  x\in (g(t),
h(t)),\, t\in (0,T]\}$;
\item[{\rm (ii)}] $0\leq V_i\leq N^*_v$  and $0\leq H_i\leq N^*_h$ for $g(t)\leq x\leq h(t),\, 0< t\leq T_0$;
 \item[{\rm (iii)}] There exists a positive constant $C_1$ such that $ 0<-g'(t),\, h'(t)\leq C_1$ for $t>0$.

\end{itemize}
\end{thm}
\bpf We only sketch the proof here since these results are basic and the proof is standard.
The local existence and uniqueness of the solution come from the contraction mapping theorem, see the detailed arguments in \cite{CF} or \cite{DL}.

Let $(V_i, H_i; g, h)$ be a solution to problem \eqref{a3} defined for $t\in [0,T_0]$ for some $T_0>0$.
It is easy to see that $0\leq V_i\leq N^*_v$  and $0\leq H_i\leq N^*_h$ for $g(t)\leq x\leq h(t),\, 0< t\leq T_0$. Noting that $H_{i0}>0$ in $(-h_0, h_0)$ and using the strong maximum principle to the equations of (\ref{a3}) in $\{(x,t):\, g(t)<x<h(t),\, 0\leq t\leq T_0\}$, we immediately obtain
 $$H_i(x, t)>0\;\; \textrm{for} \ g(t)< x<h(t),\, 0< t\leq T_0.$$
Moreover, by Hopf's lemma, we have
 $$\frac {\partial H_i}{\partial x}(g(t), t)>0\; \textrm{and} \; \frac {\partial H_i}{\partial x}(h(t), t)<0\;\; \textrm{for} \; 0<t\leq T_0.$$
 Hence $g'(t)<0$ and $h'(t)>0$ for $t\in (0, T_0]$ by using the free boundary conditions in (\ref{a3}).
Furthermore, comparison principle can be used to show that $g'(t)\geq -C_1$ and $h'(t)\leq C_1$ for $t\in (0, T_0]$ and
some $C_1$. The proof is similar to that of Lemma 2.2 in \cite{DL} with $C_1=2MN_h^*\mu D_h$ and
$$M=\max\left\{ \sqrt{\frac{\beta_h N^*_vN^*_h}{2D_h}},\
\frac{4\|H_{i0}\|_{C^1([-h_0,h_0])}}{3N^*_h}\right\}.$$

The global existence is from the uniqueness of the local solution, Zorn's lemma and the uniform estimates of $||V_i||_{C(D_{T_0})}$, $||H_i||_{C(D_{T_0})}$ and $||g'(t)||_{C((0, T_0])}, ||h'(t)||_{C((0,T_0])}$, the above estimates on $D_{T_0}$ are independent of $T_0$, therefore, these arguments hold for any $t>0$.
 \epf

Theorem \ref{exist} (iii) shows that the left free boundary for problem (\ref{a3}) is strictly monotone decreasing and the right is increasing. Ecologically, it means that the infection area which contains infected birds is always gradually expanding.

Now we present a comparison principle for problem (\ref{a3}), which can be used to estimate $V_i, H_i$ and the free
boundaries $x=g(t), x=h(t)$. We first give the following definition.

\begin{defi} The vector $(\overline V_i, \overline H_i; \overline g, \overline h)$ in
$[C^{2,1}({D}_{1T})\times
C(\overline D_{1T}))]^2\times [C^{1}([0,T])]^2$ is called an upper solution of problem \eqref{a3} if
$0\leq \overline V_i \leq N^*_v$, $0\leq \overline H_i\leq N^*_h]$ and
\begin{eqnarray*}
\left\{
\begin{array}{lll}
\dfrac{\partial \overline V_i}{\partial t}-D_v \dfrac{\partial^2 \overline V_i}{\partial x^2}\geq
\dfrac{\beta_v(N_v^*-\overline V_i)\overline H_i}{N_h^*} -r_v(1-q) \overline V_i,\; & \overline g(t)<x<\overline h(t),\, 0<t\leq T, \\
\dfrac{\partial \overline H_i}{\partial t}-D_h \dfrac{\partial^2\overline H_i}{\partial x^2}\geq
\dfrac{\beta_h \overline V_i (N^*_h-\overline H_i)}{N_h^*}-(d_h+\gamma_h)\overline H_i,\; & \overline g(t)<x<\overline h(t), \, 0<t\leq T,\\
\overline V_i(x,t)=\overline H_i(x, t)=0,&x=\overline g(t)\, \textrm{or}\, x=\overline h(t),\, 0<t\leq T,\\
\overline g(0)\leq -h_0,\; \overline g'(t)\leq -\mu D_h \frac{\partial \overline H_i}{\partial x}(\overline g(t), t), & 0<t\leq T, \\
\overline h(0)\geq h_0, \; \overline h'(t)\geq -\mu D_h \frac{\partial \overline H_i}{\partial x}(\overline h(t), t), & 0<t\leq T,\\
\overline V_i(x,0)\geq V_{i0}(x),\; \overline H_i(x,0)\geq H_{i0}(x), &-h_0\leq x\leq h_0.
\end{array} \right.
\end{eqnarray*}
$(\underline V_i, \underline H_i; \underline g, \underline h)$
in $[C^{2,1}({D}_{2T})\times
C(\overline D_{2T}))]^2\times [C^{1}([0,T])]^2$ is a lower solution if
 all the  inequalities in the obvious places are reverse, where $D_{1T}=\{(x,t): \overline g(t)<x<\overline h(t)),\, 0<t\leq T\}$
and $D_{2T}=\{(x,t): \underline g(t)<x<\underline h(t)),\, 0<t\leq T\}$.
\end{defi}

In what follows, we shall exhibit the comparison principle, and the proof is similar to that of Lemma 3.5 in \cite{DL, DL3}.
\begin{lem} (The Comparison Principle) Let $(\overline V_i, \overline H_i; \overline g, \overline h)$ and
$(\underline V_i, \underline H_i; \underline g, \underline h)$ be the upper and lower solutions of problem $(\ref{a3})$.
 Then the solution $(V_i, H_i,; g, h)$ satisfies
$$\overline g(t)\leq g(t)\leq \underline g(t) \ \textrm{and}\ \underline h(t)\leq h(t)\leq\overline h(t)\, \textrm{for}\ 0<t\leq T,$$
$$V_i\leq
\overline V_i\ \textrm{and}\ H_i\leq \overline H_i\ \textrm{for}\ g(t)\leq x\leq h(t),\, 0<t\leq T,$$
$$\underline V_i\leq
V_i\ \textrm{and}\ \underline H_i\leq H_i\ \textrm{for}\ \underline g(t)\leq x\leq \underline h(t),\, 0<t\leq T.$$
\end{lem}

\section{Basic reproduction numbers}
In this section, we present the basic reproduction numbers and their properties
for the simplified version of the WNv model subject to different environmental (boundary) settings.

If the environment is homogeneous, the governing system is given by
\begin{equation}\label{a31}
\left \{
\begin{array}{lll}
\dfrac{\textrm{d} V_i}{\textrm{d} t}=\dfrac{\beta_v(N_v^*-V_i)H_i}{N_h^*} -r_v(1-q) V_i,\; &  t>0, \\
\dfrac{\textrm{d} H_i}{\textrm{d} t}=\dfrac{\beta_h V_i (N^*_h-H_i)}{N_h^*}-(d_h+\gamma_h) H_i,\; & t>0.
\end{array} \right.
\end{equation}
It follows from a direct calculation (\cite{VW}) that the usual basic reproduction number $R_0$ for problem (\ref{a31}) is determined by
$$R_0=\ \sqrt{\frac {\beta_v\beta_hN_v^*}{r_v(1-q)N_h^*(d_h+\gamma_h)}}\, .$$

It is not difficult to
see that there is the only trivial steady state ( disease-free equilibrium ) $(0, 0)$
if $R_0\leq 1$, while if $R_0>1$, there exists a unique positive
steady state ( endemic equilibrium ) $(V^*_i, H^*_i)$ with
$$V_i^*=\frac {\beta_v\beta_hN_v^*-r_v(1-q)N_h^*(d_h+\gamma_h)}{\beta_v\beta_h+r_v(1-q)\beta_h},\ H_i^*=\frac{\beta_hV_i^*}{\frac{\beta_hV_i^*}{N_h^*}+d_h+\gamma_h}.$$
Moreover, we can prove that $(0, 0)$ is globally asymptotically stable if $R_0\leq 1$ by using Lyapunov functional, and the positive
steady state $(V^*_i, H^*_i)$ is locally asymptotically stable if $R_0>1$ by using standard linearization and spectral analysis, as illustrated in Fig. 2.
The positive steady state is also globally asymptotically stable by Poincar$\acute{e}$-Bendixson theorem, see Proposition 2.1
in \cite{Lewis05}, where the system with different coefficients is studied.

 \begin{figure}[htbp] \centering
\includegraphics[width=9cm]{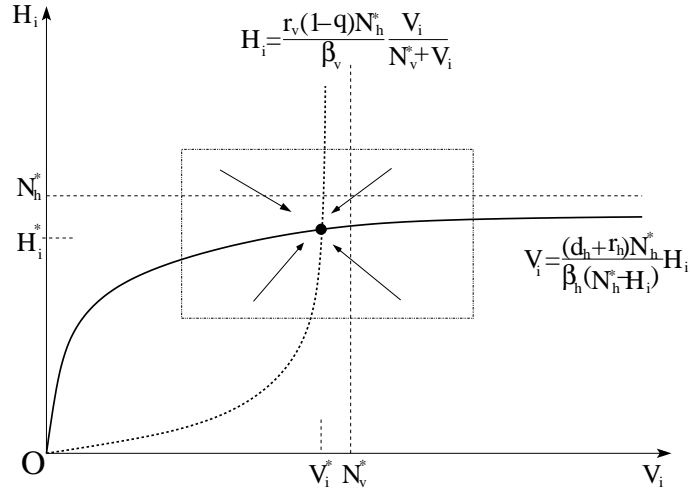}
   \caption{Phase portrait for the governing system \eqref{a31}  in $(V_i, H_i)$ plane. }\label{t4}
\end{figure}

If we consider the transmission of WNv in a bounded region, denoted by $\Omega$,
 then there are two special cases on the fixed boundary:
one case is that both vector mosquitoes and host birds do not cross the boundary,
the other case is that there are neither infected mosquitoes nor infected birds on the boundary of $\Omega$.

In a bounded region $\Omega$ with $\partial \Omega \in C^1$,
if there are neither vector mosquitos nor host birds crossing the boundary, the
corresponding spatially-dependent WNv transmission model
can be written as
\begin{equation}\label{a32}
\left \{
\begin{array}{lll}
\dfrac{\partial V_i}{\partial t}-D_v \Delta V_i=\dfrac{\beta_v(N_v^*-V_i)H_i}{N_h^*} -r_v(1-q) V_i,\; & x\in \Omega,\, t>0, \\
\dfrac{\partial H_i}{\partial t}-D_h \Delta H_i=\dfrac{\beta_h V_i (N^*_h-H_i)}{N_h^*}-(d_h+\gamma_h) H_i,\; & x\in \Omega, \, t>0,\\
\frac {\partial V_i}{\partial \eta}=\frac {\partial H_i}{\partial \eta}=0,\, &  x\in \partial \Omega, \ t>0.
\end{array} \right.
\end{equation}

Let $R_0^N=R_0^N(\Omega)$ be the positive principal eigenvalue to the problem
\begin{eqnarray}\label{a321}
\left\{
\begin{array}{lll}
-D_v \Delta \phi=\frac{\beta_v N_v^*}{N_h^*R_0^N}\psi -r_v(1-q)\phi,\; & x\in \Omega,\\
-D_h\Delta \psi= \frac{\beta_h}{R_0^N}\phi-(d_h+\gamma_h)\psi,\; & x\in \Omega,  \\
\frac{\partial\phi}{\partial \eta}(x)=\frac{\partial\psi}{\partial \eta}(x)=0, &x\in \partial \Omega.
\end{array} \right.
\end{eqnarray}
The principal eigenvalue is the only positive eigenvalue admitting a unique positive eigenfunction $(\phi,\psi)$ (subject to a constant multiple).

To get the existence of $R_0^N(\Omega)$, we consider the following eigenvalue problem,
\begin{eqnarray}\label{a322}
\left\{
\begin{array}{lll}
-D_v \Delta \phi=\frac{\beta_v N_v^*}{N_h^*R}\psi -r_v(1-q)\phi+\mu\phi,\; & x\in \Omega,\\
-D_h\Delta \psi= \frac{\beta_h}{R}\phi-(d_h+\gamma_h)\psi+\mu\psi,\; & x\in \Omega,  \\
\frac{\partial\phi}{\partial \eta}(x)=\frac{\partial\psi}{\partial \eta}(x)=0, &x\in \partial \Omega.
\end{array} \right.
\end{eqnarray}
Recalling that the system is strongly cooperative, one may argue as in \cite{ALG, LG} to show that, for any fixed positive $R$, there is a unique value $\mu=\mu_1(R)$ such that \eqref{a322} admits a unique positive
solution $(\phi,\psi)$ (subject to a constant multiple). Such a value $\mu$ is known as the principal eigenvalue of \eqref{a322}. Moreover, $\mu_1(R)$ is
continuous and strictly increasing. The existence of $R_0^N(\Omega)$ follows from the fact that $\lim_{R\to 0^+}\mu_1(R)<0$ and $\lim_{R\to +\infty}\mu_1(R)>0$.
Since all coefficients in \eqref{a322} are constants, it is easy to check that $R_0^N(\Omega)=R_0$.

For the second case with bounded $\Omega$,
if there is neither infected vector mosquitos nor host birds on the boundaries,
then the corresponding spatially-dependent WNv transmission model becomes
\begin{equation}\label{a33}
\left \{
\begin{array}{lll}
\dfrac{\partial V_i}{\partial t}-D_v \Delta V_i=\dfrac{\beta_v(N_v^*-V_i)H_i}{N_h^*} -r_v(1-q) V_i,\; & x\in \Omega,\, t>0, \\
\dfrac{\partial H_i}{\partial t}-D_h \Delta H_i=\dfrac{\beta_h V_i (N^*_h-H_i)}{N_h^*}-(d_h+\gamma_h) H_i,\; & x\in \Omega, \, t>0,\\
V_i(x,t)= H_i(x,t)=0,\, &  x\in \partial \Omega, \ t>0.
\end{array} \right.
\end{equation}

We now introduce the basic reproduction number $R_0^D=R_0^D(\Omega)$ for model \eqref{a33} by the positive principal eigenvalue to the problem
\begin{eqnarray}\label{a32-1}
\left\{
\begin{array}{lll}
-D_v \Delta \phi=\frac{\beta_v N_v^*}{N_h^*R_0^D}\psi -r_v(1-q)\phi,\; & x\in \Omega,\\
-D_h\Delta \psi= \frac{\beta_h}{R_0^D}\phi-(d_h+\gamma_h)\psi,\; & x\in \Omega,  \\
\phi(x)=\psi(x)=0, &x\in \partial \Omega,
\end{array} \right.
\end{eqnarray}
where the corresponding eigenfunction $(\phi,\psi)$ is unique, up to a positive multiplicative constant, and $\phi(x)>0,\psi(x)>0$ for all $x\in \Omega$.

Normally, we cannot use variational methods to treat eigenvalue problems for coupled systems, though variational methods
are proved to be effective for most scalar eigenvalue problems.
Thanks to the assumption that all coefficients are constant, we can provide an explicit formula for $R_0^D$ and a similar result to Lemma 2.3 in \cite{HH}.
\begin{lem} Problem \eqref{a32-1} admits a unique positive principal eigenvalue $R_0^D(\Omega)$ determined by
\begin{equation}\label{basic2}
R_0^D(\Omega)=\sqrt{\frac{ \frac {\beta_v\beta_hN_v^*}{N_h^*}}{ [D_h\lambda^*+(d_h+\gamma_h)][D_v\lambda^*+r_v(1-q)]}}\ ,
\end{equation}
where $\lambda^*$ is the principal eigenvalue of $-\Delta$ in $\Omega$ with null Dirichlet boundary condition.
Moreover, there exist $\lambda_0$ and $\delta_0>0$ such that
$1-R_0^D$ has the same sign as $\lambda_0$, and $(\phi,\psi):=(\delta_0 \psi^*, \psi^*)$ satisfies
\begin{eqnarray}
\left\{
\begin{array}{lll}
-D_v \Delta \phi=\frac{\beta_v N_v^*}{N_h^*}\psi -r_v(1-q)\phi+\lambda_0 \phi,\; & x\in \Omega,\\
-D_h\Delta \psi= \beta_h\phi-(d_h+\gamma_h)\psi+\lambda_0 \psi,\; & x\in \Omega,  \\
\phi(x)=\psi(x)=0, &x\in \partial \Omega,
\end{array} \right.
\label{B1f}
\end{eqnarray}
where $\psi^*$ is the eigenfunction corresponding to the principal eigenvalue ($\lambda^*$) of $-\Delta$ in $\Omega$ with null Dirichlet boundary condition.
\end{lem}
\bpf
Let
$$R^*=\frac{ \frac {\beta_v\beta_hN_v^*}{N_h^*}}{ [D_h\lambda^*+(d_h+\gamma_h)][D_v\lambda^*+r_v(1-q)]},$$
$$\phi^*=\frac{\beta_vN_v^*}{N_h^*\sqrt{R^*}[D_v\lambda^*+r_v(1-q)]}\psi^*.$$
Then we know that $(\phi^*,\psi^*)$ is a positive solution of problem (\ref{a32-1}) with $R_0^D=\sqrt{R^*}$,
and \eqref{basic2} follows directly from the uniqueness of the principal eigenvalue of (\ref{a32-1}).

It is easy to verify that $(\phi,\psi):=(\delta_0 \psi^*, \psi^*)$ satisfies \eqref{B1f} if
\begin{eqnarray}
\left\{
\begin{array}{lll}
D_v \lambda^* \delta_0=\frac{\beta_v N_v^*}{N_h^*} -r_v(1-q)\delta_0+\lambda_0 \delta_0,\; & \\
D_h\lambda^*= \beta_h\delta_0-(d_h+\gamma_h)+\lambda_0.\; &
\end{array} \right.
\label{B1af}
\end{eqnarray}
Denote
$$A=D_v \lambda^*+r_v(1-q),\ \ B=D_h\lambda^*+(d_h+\gamma_h)$$
for simplicity, and choose
$$\lambda_0=\frac{(A+B)-\sqrt{(A+B)^2+4AB(R^*-1)}}2,\ \ \delta_0=\frac {B-\lambda_0}{\beta_h},$$
then \eqref{B1af} holds. Moreover, $\textrm{sign}(1-R_0^D)=\textrm{sign}(1-R^*)=\textrm{sign} \lambda_0$ and
$$\delta_0=\frac {2B-(A+B)+\sqrt{(A+B)^2+4AB(R^*-1)}}{2\beta_h}>\frac {2B-(A+B)+|A-B|}{2\beta_h}\geq 0.$$
\epf
\bigskip

A straightforward calculation shows that
\begin{lem} The following assertions hold.
\begin{enumerate}
\item[$(a)$]  $R_0^D$ is a positive and monotonically decreasing function of $D_h$ and $D_v$; $R_0^D\to R_0$ as $D_h\to 0$ and $D_v\to 0$,
 $R_0^D\to 0$ as $D_h\to \infty$ or $D_v\to \infty$;

\item[$(b)$] Let $B_h$ be a ball with the radius $h$. $R_0^D(B_{h})$ is strictly monotonicaly increasing function of $h$, that is if $h_1<h_2$, then $R_0(B_{h_1})<R_0(B_{h_2})$.
Moreover, $R_0^D(B_{h})\to R_0$ as $h\to \infty$ and $R_0^D(B_{h})\to 0$ as $h\to 0$;

\item [$(c)$] If $\Omega =(-h_0, h_0)$, then
$$R_0^D(\Omega)=\sqrt{\frac{ \frac {\beta_v\beta_hN_v^*}{N_h^*}}{ [D_h(\frac \pi{2h_0})^2+(d_h+\gamma_h)][D_v(\frac \pi{2h_0})^2+r_v(1-q)]}}\, .$$
\end{enumerate}
\end{lem}

Note that the domain $(g(t), h(t))$ for the free boundary problem \eqref{a3} is changing with $t$, so the corresponding threshold value is not a constant.
Now we introduce the threshold value $R_0^F(t)$ for the free boundary problem \eqref{a3} by
$$R_0^F(t)=R_0^D((g(t),h(t))).$$
Since $R_0^F(t)$ involves regional characteristic and time, we then call it {\bf the spatial-temporal risk index}, instead of the basic reproduce number.
Recalling that all coefficients in \eqref{a3} are constants, we then have
$$R_0^F(t)=\sqrt{\frac{ \frac {\beta_v\beta_hN_v^*}{N_h^*}}{ [D_h(\frac \pi{h(t)-g(t)})^2+(d_h+\gamma_h)][D_v(\frac \pi{h(t)-g(t)})^2+r_v(1-q)]}}.
$$
It follows from Theorem 3.1 $(iii)$ and Lemma 4.2 that
\begin{lem} $R_0^F(t)$ is strictly monotonically increasing function of $t$, that is if $t_1<t_2$, then $R_0^F(t_1)<R_0^F(t_2)$.
Moreover, if $h(t)-g(t)\to \infty$ as $t\to \infty$, then $R_0^F(t)\to R_0$ as $t\to \infty$.
\end{lem}

\section{The scenario of vanishing}

In this section, we will consider the vanishing scenario of the WNv. First we know that the two free boundary fronts $g(t)$ and $h(t)$
have the monotonicity described in Theorem \ref{exist}, the next theorem shows that the two free boundary
fronts $x=g(t)$ and $x=h(t)$ are both finite or infinite simultaneously.

\begin{thm} Let $(V_i, H_i; g, h)$ be a solution to problem $(\ref{a3})$
defined for $t\in[0, +\infty)$ and $x\in[g(t), h(t)]$. Then for $t\in[0,
+\infty)$ we have
$$-2h_0<g(t)+h(t)<2h_0.$$
\end{thm}
\bpf By continuity we know that $g(t)+h(t)>-2h_0$ holds for small $t>0$.
Let
$$T:=\sup\{s: g(t)+h(t)>-2h_0 \mbox{ for all }  t\in[0,s)\}.$$
As in \cite{DLou}, we claim
that $T=\infty$. Otherwise, we have $0<T<\infty$  and
$$g(t)+h(t)>-2h_0 \mbox{ for } t\in[0,T),\ g(T)+h(T)=-2h_0.$$
Hence
\begin{eqnarray}
g'(T)+h'(T)\leq0. \label{Hq}
\end{eqnarray}

To reach a contradiction, we consider two functions
$$\begin{array}{rcl}
w(x, t):&=&V_i(x,t)-V_i(-x-2h_0, t),\\
z(x, t):&=&H_i(x, t)-H_i(-x-2h_0, t)
\end{array}
$$
over the region $$\Lambda:=\{(x, t):\ x\in[g(t), -h_0],\, t\in[0, T]\}.$$
It is not difficult to verify that
the pair $(w, z)$ is
well-defined for $(x, t)\in\Lambda$ since $-h_0\leq -x-2h_0\leq -g(t)-2h_0\leq h(t)$, and the pair satisfies
$$w_t-D_v w_{xx}=c_{11}(x,t)w+c_{12}(x,t)z, \ \ \ \ \ \mbox{ for } g(t)<x<-h_0,\ 0<t\leq T,$$
$$z_t-D_h z_{xx}=c_{21}(x,t)w+c_{22}(x,t)z,  \ \ \ \ \ \ \mbox{ for } g(t)<x<-h_0,\ 0<t\leq T$$
with some $c_{12}\geq 0, c_{21}\geq 0$ and $c_{ij}\in L^\infty(\Lambda)$ for $i, j=1, 2$, and
$$w(-h_0, t)=z(-h_0,t)=0,\, w(g(t), t)\leq 0,\, z(g(t), t)<0 \ \ \ \ \ \mbox{ for } 0<t<T.$$ Moreover,
$$
\begin{array}{rcl}
z(g(T), T)&=&H_i(g(T), T)-H_i(-g(T)-2h_0, T)\\
         &=&H_i(g(T), T)-H_i(h(T), T)=0.
         \end{array}
         $$
Applying the proof for the strong maximum principle and the Hopf's lemma, we deduce
$$w(x,t)<0,\, z(x,t)<0  \mbox{ in } (g(t),-h_0)\times (0, T] \mbox{ and } z_x(g(T), T)<0.$$
But we know
$$z_x(g(T), T)=\frac {\partial H_i}{\partial x}(g(T), T)+\frac {\partial H_i}{\partial x}(
h(T), T)=-[g'(T)+h'(T)]/(D_h\mu),$$ which implies $$g'(T)+h'(T)>0.$$
Therefore there is a
contradiction to (\ref{Hq}). Hence we prove $$g(t)+h(t)>-2h_0  \ \ \  \mbox{
for all } t>0.$$

Analogously we can prove $g(t)+h(t)<2h_0 \mbox{ for all } t>0$ by considering
$$W(x, t):=V_i(x,t)-V_i(2h_0-x, t),\ Z(x, t):=H_i(x, t)-H_i(2h_0-x, t)$$ over the
region $\Lambda':=[h_0, h(t)]\times [0, T']$ with
$T':=\sup\{s: g(t)+h(t)<2h_0 \mbox{ for all }  t\in[0,s)\}$. The
proof is completed.  \epf

It follows from Theorem \ref{exist} that  $x=g(t)$ is monotonically decreasing and $x=h(t)$ is monotonically increasing, therefore there exist $\displaystyle g_\infty\in [-\infty, 0)$ and $h_\infty\in (0, +\infty]$ such that
$\displaystyle \lim_{t\to +\infty} \ g(t)=g_\infty$ and $\lim_{t\to +\infty} \ h(t)=h_\infty$. Next we discuss the properties of the free boundary.
Since the transmission of the virus
depends on whether or not $h_\infty-g_\infty=\infty$ and $\displaystyle\lim_{t\to
+\infty} \ ||H_i(\cdot, t)||_{C([g(t), h(t)])}=0$, we give the following definitions representing two different scenarios of the virus transmission:

\begin{defi}
The virus is {\bf vanishing} if $h_\infty-g_\infty <\infty$ and
 $\displaystyle \lim_{t\to +\infty} \ ||H_i(\cdot, t)||_{C([g(t),h(t)])}=0$; the virus is {\bf spreading} if $h_\infty -g_\infty=\infty$ and
 $\displaystyle \limsup_{t\to +\infty} \ ||H_i(\cdot, t)||_{C([g(t),h(t)])}>0$.
\end{defi}

The next result shows that if $h_\infty-g_\infty <\infty$, then vanishing scenario will happen.
\begin{lem}  If $h_\infty-g_\infty <\infty$, then we have
$$\lim_{t\to
+\infty} \ ||H_i(\cdot, t)||_{C([g(t), h(t)])}=\lim_{t\to
+\infty} \ ||V_i(\cdot, t)||_{C([g(t), h(t)])}=0.$$
\label{vanishlm}
\end{lem}
\bpf
Assume  that $\displaystyle \limsup_{t\to
+\infty} \ ||H_i(\cdot, t)||_{C([g(t), h(t)])}=\delta>0$ by contradiction. Then there exists a sequence $(x_k, t_k )$
in $(g(t), h(t))\times (0, \infty)$
such that $H_i(x_k,t_k)\geq \frac{\delta}{2}$ for all $k \in \mathbb{N}$, and $t_k\to \infty$ as $k\to \infty$.
Since $-\infty<g_\infty<g(t)<x_k<h(t)<h_\infty<\infty$, we then obtain that there exists a subsequence of $\{x_n\}$ converging
to $x_0\in (g_\infty, h_\infty)$. Without loss of generality, we assume $x_k\to x_0$ as $k\to \infty$.

Let $W_k(x,t)=V_i(x,t_k+t)$ and $Z_k(x,t)=H_i(x,t_k+t)$  for
$x\in (g(t_k+t), h(t_k+t)), t\in (-t_k, \infty)$. As in \cite{FS},
it follows from the parabolic regularity that  $\{(W_k, Z_k)\}$ has a subsequence $\{(W_{k_i}, Z_{k_i})\}$ such that
$(W_{k_i}, Z_{k_i})\to (\tilde W, \tilde Z)$ as $i\to \infty$ and $(\tilde W, \tilde Z)$ satisfies
\begin{eqnarray*} \left\{
\begin{array}{lll}
\tilde W_t-D_v \tilde W_{xx}=\dfrac{\beta_v(N_v^*-\tilde W)\tilde Z}{N_h^*} -r_v(1-q) \tilde W,\; & g_\infty<x<h_\infty,\ t\in (-\infty, \infty),  \\
\tilde Z_t-D_h \tilde Z_{xx}=\dfrac{\beta_h \tilde W (N^*_h-\tilde Z)}{N_h^*}-(d_h+\gamma_h)\tilde Z,\; &\ g_\infty<x<h_\infty, \ t\in (-\infty, \infty).
\end{array} \right.
\end{eqnarray*}
Note that $\tilde Z(x_0, 0)\geq \delta/2$, therefore $\tilde Z>0$ in $ (g_\infty, h_\infty)\times(-\infty, \infty)$.

 Using the similar method to prove Hopf's lemma at the point $(h_\infty, 0)$ yields
that $\tilde Z_x(h_\infty, 0 )\leq -\sigma_0$ for some $\sigma_0>0$.

On the other hand, since $-g(t)$ and $h(t)$ are increasing and bounded, it follows from standard $L^p$ theory and the Sobolev imbedding
theorem (\cite{LSU}) that, for any $0<\alpha <1$, there exists a constant $\tilde C$
depending on $\alpha$, $h_0$, $\|V_{i0}\|_{C^{2}[-h_0, h_0]}$, $\|H_{i0}\|_{C^{2}[-h_0, h_0]}$ and $g_\infty, h_\infty$ such that
\begin{eqnarray}
\|H_i\|_{C^{1+\alpha,
(1+\alpha)/2}([-g(t), h(t)]\times[\tau, \tau+1])}\leq \tilde C\label{Bg1}
\end{eqnarray}
for any $\tau\geq 1$. Noting that $\tilde C$ is independent of $\tau$ and using the free boundary conditions in (\ref{a3}), we then achieve
\begin{eqnarray}
\|H_i(\cdot, t)\|_{C^{1}([g(t), h(t)])}\leq \hat C,\ t\geq 1 ,\label{est-2}\\
||h'||_{C^{\alpha/2}([1, +\infty))},\, ||g'||_{C^{\alpha/2}([1, +\infty))}\leq \hat C.
\label{est-1}
\end{eqnarray}

Now, since $\|h'\|_{C^{\alpha/2}([1,\infty))}\leq
\hat C$ and $h(t)$ is bounded,  we then have $h'(t)\to 0$ as $t\to \infty$, that is,
$\frac {\partial H_i}{\partial x}(h(t_k),t_k)\to 0$ as $t_k\to \infty$ by the free boundary condition. Moreover, using (\ref{est-2}) gives that
$$
\frac {\partial H_i}{\partial x}(h(t_k),t_k+0)=(Z_k)_x(h(t_k),0)\to \tilde Z_x(h_\infty,0),  \ \ \ \ \ \mbox{as} \ \ k\to \infty $$
which leads to a contradiction to the fact that  $\tilde Z_x(h_\infty,0)\leq -\sigma_0<0$.
 Thus  we have
 $$\displaystyle \lim_{t\to +\infty} \ ||H_i(\cdot,t)||_{C([g(t),h(t)])}=0.$$

 The above limit implies that for any $\varepsilon>0$, there exists a $T>0$ such that $0\leq H_i(x,t)\leq \varepsilon$
 for $x\in [g(t), h(t)]$ and $t\geq T$. Note that $V_i$ satisfies
 $$\frac{\partial V_i}{\partial t}-D_v \frac{\partial^2 V_i}{\partial x^2}\leq \frac{\beta_v N_v^*}{N_h^*}\varepsilon -r_v(1-q) V_i,\; g(t)<x<h(t),\, \ \  t\in [T,\infty).$$
 Therefore $\displaystyle \limsup_{t\to +\infty} \ ||V_i(\cdot,t)||_{C([g(t),h(t)])}\leq \frac {\beta_v N_v^*}{N_h^*r_v(1-q)}\varepsilon$.
Since $\varepsilon$ is sufficiently small, we have
$\displaystyle \lim_{t\to +\infty} \ ||V_i(\cdot,t)||_{C([g(t),h(t)])}=0$.
\epf

\bigskip
Next we exhibit sufficient conditions, under which the transmission of the virus is in a vanishing scenario.

\begin{thm} If $R_0\leq 1$, then $h_\infty-g_\infty<\infty$ and vanishing happens.
\label{vanish}
\end{thm}
\bpf
By Lemma \ref{vanishlm}, it suffices to prove that $h_\infty -g_\infty<+\infty$.
Direct calculations yield
\begin{eqnarray*}& &\frac{\textrm{d}}{\textrm{d} t}\int_{g(t)}^{h(t)}[V_i+\frac {r_v(1-q)}{\beta_h}H_i](x,t)\textrm{d}x\\[1mm]
&=&\int_{g(t)}^{h(t)}[\frac{\partial V_i}{\partial t}+\frac {r_v(1-q)}{\beta_h}\frac{\partial H_i}{\partial t}](x, t)\textrm{d}x+h'(t)[V_i+\frac {r_v(1-q)}{\beta_h}H_i](h(t), t)\\
& &-g'(t)[V_i+\frac {r_v(1-q)}{\beta_h}H_i](g(t), t)\\
&=&\int_{g(t)}^{h(t)}(D_v\frac{\partial^2 V_i}{\partial x^2}+\frac {r_v(1-q)}{\beta_h} D_h\frac{\partial^2 H_i}{\partial x^2})\textrm{d}x\\
& &+\int_{g(t)}^{h(t)}H_i\big(\frac{\beta_vN^*_v}{N^*_h}-\frac{r_v(1-q)}{\beta_h}(d_h+\gamma_h)-
\frac{(\beta_v+r_v(1-q))}{N^*_h}V_i\big)\textrm{d}x\\
&=&-\frac {r_v(1-q)}{\mu \beta_h}(h'(t)-g'(t))+D_v\big(\frac{\partial V_i}{\partial x}(h(t),t)-\frac{\partial V_i}{\partial x}(g(t),t)\big)\\
& &+\int_{g(t)}^{h(t)}H_i\big[\frac{\beta_vN^*_v}{N^*_h}(1-\frac{1}{R_0})-
\frac{(\beta_v+r_v(1-q))}{N^*_h}V_i\big]\textrm{d}x.
\end{eqnarray*}
Recalling that $R_0\leq 1$, $\frac{\partial V_i}{\partial x}(h(t),t)\leq 0$ and $\frac{\partial V_i}{\partial x}(g(t),t)\geq 0$, and integrating from $0$ to $t\,(>0)$ give
$$\begin{array}{rl}
 \int_{g(t)}^{h(t)}[V_i+\frac {r_v(1-q)}{\beta_h}H_i](x, t)\textrm{d}x
  &\leq \int ^{h(0)}_{g(0)}[V_i+\frac {r_v(1-q)}{\beta_h}H_i](x, 0)\textrm{d}x  \\[2ex]
  &\quad +\frac {r_v(1-q)}{\mu \beta_h}(h(0)-g(0))-\frac {r_v(1-q)}{\mu \beta_h}(h(t)-g(t)), \quad t\geq 0.
\end{array}
$$
We then have
$$\frac {r_v(1-q)}{\mu \beta_h}(h(t)-g(t)) \leq \int ^{h(0)}_{g(0)}[V_i+\frac {r_v(1-q)}{\beta_h}H_i](x, 0)\textrm{d}x
+\frac {r_v(1-q)}{\mu \beta_h}(h(0)-g(0))$$
for $t\geq 0$, which in turn gives that $h_\infty-g_\infty<\infty$. Therefore, the virus is vanishing.
\epf

\begin{thm} Suppose $R_0^F(0)(:=R_0^D((-h_0, h_0)))<1$.  Then $h_\infty-g_\infty<\infty$ and
$$\lim_{t\to
+\infty} \ ||H_i(\cdot, t)||_{C([g(t), h(t)])}=\lim_{t\to
+\infty} \ ||V_i(\cdot, t)||_{C([g(t), h(t)])}=0$$
 if $||V_{i0}(x)||_{C([-h_0, h_0])}$
and $||H_{i0}(x)||_{C([-h_0, h_0])}$ are sufficiently small.
\end{thm}
\bpf We are going to construct a suitable upper solution for problem (\ref{a3}).
Note that $R_0^D((-h_0, h_0))<1$, it follows from Lemma 4.1 that there exist $\lambda_0>0$, $\delta_0>0$ and $\psi(x)>0$ in $(-h_0, h_0)$ such that
\begin{eqnarray}
\left\{
\begin{array}{lll}
-D_v \Delta \phi=\frac{\beta_v N_v^*}{N_h^*}\psi -r_v(1-q)\phi+\lambda_0 \phi,\; & -h_0<x<h_0,\\
-D_h\Delta \psi= \beta_h\phi-(d_h+\gamma_h)\psi+\lambda_0 \psi,\; &-h_0<x<h_0,  \\
\phi(x)=\psi(x)=0, &x=\pm h_0.
\end{array} \right.
\label{B1f1}
\end{eqnarray}
where $\phi(x)=\delta_0 \psi(x)$. Therefore, there exists a small $\delta >0$ such that
$$-\delta +(\frac 1{(1+\delta)^2}-1)[\dfrac{\beta_vN_v^*}{N_h^*\delta_0} -r_v(1-q)]+\frac 1{(1+\delta)^2}\lambda_0\geq 0$$
and
$$-\delta +(\frac 1{(1+\delta)^2}-1)[\beta_h\delta_0-(d_h+\gamma_h)]+\frac 1{(1+\delta)^2}\lambda_0\geq 0.$$

As in \cite{DL}, if we set
$$\sigma (t)=h_0(1+\delta-\frac \delta 2 e^{-\delta t}), \  t\geq 0,$$
and
$$\overline V(x, t)=\varepsilon e^{-\delta t}\phi(xh_0/\sigma (t)), \ -\sigma(t)\leq
x\leq \sigma(t),\ t\geq 0,$$
$$\overline H(x, t)=\varepsilon e^{-\delta t}\psi(xh_0/\sigma (t)), \ -\sigma(t)\leq
x\leq \sigma(t),\ t\geq 0,$$
then straightforward computations lead to the following
\begin{eqnarray*}
& &\overline V_t-D_v \dfrac{\partial^2 \overline V}{\partial x^2}-
\dfrac{\beta_v(N_v^*-\overline V)\overline H}{N_h^*} +r_v(1-q) \overline V\\
& &=-\delta \overline V-\varepsilon e^{-\delta t}\phi'\frac{xh_0\sigma'(t)}{\sigma^2(t)}+(\frac{h_0}{\sigma(t)})^2
[\dfrac{\beta_vN_v^*}{N_h^*\delta_0} -r_v(1-q)+\lambda_0 ]\overline V\\
& & \quad -\dfrac{\beta_v \overline H (N^*_v-\overline V)}{N_h^*}+r_v(1-q)\overline V\\
& &\geq \overline V \{-\delta +(\frac 1{(1+\delta)^2}-1)[\dfrac{\beta_vN_v^*}{N_h^*\delta_0} -r_v(1-q)]+\frac 1{(1+\delta)^2}\lambda_0\}\geq 0,
\end{eqnarray*}
\begin{eqnarray*}
& &\overline H_t-D_h \dfrac{\partial^2\overline H}{\partial x^2}-
\dfrac{\beta_h \overline V (N^*_h-\overline H)}{N_h^*}+(d_h+\gamma_h)\overline H\\
& &=-\delta \overline H-\varepsilon e^{-\delta t}\psi'\frac{xh_0\sigma'(t)}{\sigma^2(t)}+(\frac{h_0}{\sigma(t)})^2[\beta_h\delta_0-(d_h+\gamma_h)+\lambda_0 ]\overline H\\
& & \quad -\dfrac{\beta_h \overline V (N^*_h-\overline H)}{N_h^*}+(d_h+\gamma_h)\overline H\\
& &\geq \overline H \{-\delta +(\frac 1{(1+\delta)^2}-1)[\beta_h\delta_0-(d_h+\gamma_h)]+\frac 1{(1+\delta)^2}\lambda_0\}\geq 0
\end{eqnarray*}
for all $t>0$ and $-\sigma (t)<x<\sigma (t)$.

On the other hand, we come to the result that
$$
\begin{array}{rcl}
&&\sigma'(t)=\displaystyle  h_0 \frac {\delta^2} 2 e^{-\delta t},  \\
&& -\overline H_x(\sigma (t),t)=\displaystyle
     -\varepsilon \frac {h_0}{\sigma (t)}\psi'(h_0)e^{-\delta t},\\
&&-\overline H_x( -\sigma (t),t)=\displaystyle
   -\varepsilon \frac {h_0}{\sigma (t)}\psi'(-h_0)e^{-\delta t}.
   \end{array}
   $$
Noticing that $\psi'(-h_0)=-\psi'(h_0)$,
we now choose $\varepsilon=-\frac {\delta^2h_0} {4\mu\psi'(h_0)}$ such that
\begin{eqnarray*}
\left\{
\begin{array}{lll}
\dfrac{\partial \overline V}{\partial t}-D_v \dfrac{\partial^2 \overline V}{\partial x^2}\geq
\dfrac{\beta_v(N_v^*-\overline V)\overline H}{N_h^*} -r_v(1-q) \overline V,\; & -\sigma(t)<x<\sigma(t),\, t>0, \\
\dfrac{\partial \overline H}{\partial t}-D_h \dfrac{\partial^2\overline H}{\partial x^2}\geq
\dfrac{\beta_h \overline V (N^*_h-\overline H)}{N_h^*}-(d_h+\gamma_h)\overline H,\; &  -\sigma(t)<x<\sigma(t), \, t>0,\\
\overline V(x,t)=\overline H(x, t)=0,&x=\pm \sigma(t)\, \, t>0,\\
-\sigma (0)<-h_0,\; -\sigma'(t)\leq -\mu D_h \frac{\partial \overline H}{\partial x}(-\sigma(t), t), & t>0, \\
\sigma(0)> h_0, \; \sigma'(t)\geq -\mu D_h \frac{\partial \overline H}{\partial x}(\sigma(t), t), & t>0.
\end{array} \right.
\end{eqnarray*}
If  $||V_{i0}||_{L^\infty}\leq \varepsilon \phi(\frac {h_0}{1+\delta/2})$
and $||H_{i0}||_{L^\infty}\leq \varepsilon \psi(\frac {h_0}{1+\delta/2})$, then for $x\in [-h_0, h_0]$,
$$V_{i0}(x)\leq \varepsilon \phi(\frac {h_0}{1+\delta/2})\leq \overline V(x, 0)$$
and
$$H_{i0}(x)\leq \varepsilon \psi(\frac {h_0}{1+\delta/2})\leq \overline H(x, 0)$$
owing to $h_0<\sigma (0)=h_0(1+\delta/2)$.
We then can apply Lemma 3.4 to conclude that $g(t)\geq -\sigma(t)$ and $h(t)\leq\sigma(t)$ for $t>0$. It
follows that $\displaystyle h_\infty-g_\infty\leq \lim_{t\to\infty}
2\sigma(t)=2h_0(1+\delta)<\infty$, and $\displaystyle \lim_{t\to
+\infty} \ ||H_i(\cdot, t)||_{C([0, h(t)])}=0$ by Lemma 5.2.
 \epf

 It follows from the above proof that we can construct a suitable upper solution so that the virus is vanishing for small $\mu$, see also Lemma 5.10 in \cite{DL}.

 \begin{thm} Suppose $R_0^F(0)(:=R_0^D((-h_0, h_0)))<1$.  Then there exists
$\overline \mu>0$ depending on $V_{i0}$ and $H_{i0}$ such that
$h_\infty-g_\infty<\infty$ and $\displaystyle \lim_{t\to
+\infty} \ ||H_i(\cdot, t)||_{C([g(t), h(t)])}=0$ when $\mu\leq \overline \mu$.
\end{thm}

\section{The scenario of spreading}
In this section, we are going to give sufficient conditions, under which the virus is in a scenario of
spreading. We first
prove that if $R_0^F(0)\geq 1$, the virus is spreading.
\begin{thm} If $R_0^F(0)(:=R_0^D((-h_0, h_0)))\geq 1$, then $h_\infty-g_\infty=\infty$ and $$\liminf_{t\to
+\infty} \ ||H_i(\cdot, t)||_{C([g(t), h(t)])}>0,$$ that is, spreading happens.
\end{thm}
\bpf We first consider the case $R_0^F(0):=R_0^D((-h_0, h_0))>1$. In this case, the following problem
\begin{eqnarray}
\left\{
\begin{array}{lll}
-D_v \Delta \phi=\frac{\beta_v N_v^*}{N_h^*}\psi -r_v(1-q)\phi+\lambda_0 \phi,\; & -h_0<x<h_0,\\
-D_h\Delta \psi= \beta_h\phi-(d_h+\gamma_h)\psi+\lambda_0 \psi,\; &-h_0<x<h_0,  \\
\phi(x)=\psi(x)=0, &x=\pm h_0.
\end{array} \right.
\label{B2f}
\end{eqnarray}
admits a positive solution $(\phi(x), \psi(x))$ with $\phi(x)=\delta_0 \psi(x)$ and $\psi(x)(=\cos (\frac{x\pi}{2h_0}))$ is the eigenfunction corresponding to the principal eigenvalue of $-\Delta$ in $(-h_0, h_0)$ with null Dirichlet boundary condition. It follows from Lemma 4.1
that $\lambda_0<0$ and $\delta_0>0$.

We are now in a position to construct a suitable lower solution to
\eqref{a3} and let
$$\underline {V}(x,t)=\delta \phi(x),\quad \underline H=\delta \psi(x)$$
for $-h_0\leq x\leq h_0$, $t\geq 0$, where $\delta$ is sufficiently small.

 Direct computations yield
\begin{eqnarray*}
& &\dfrac{\partial \underline V}{\partial t}-D_v \dfrac{\partial^2 \underline V}{\partial x^2}-
\dfrac{\beta_v(N_v^*-\underline V)\underline H}{N_h^*} +r_v(1-q) \underline V \\
& &\quad =\delta [\frac {\beta_vN_v^*}{N^*_h}\psi -r_v(1-q)\phi+\lambda_0 \phi-\frac {\beta_v(N_v^*-\delta\phi)}{N^*_h}
\psi+r_v(1-q)\phi],\\
& &\quad =\delta \phi [\lambda_0 +\frac {\delta \beta_v\psi}{N^*_h}],\\
& &\dfrac{\partial \underline H}{\partial t}-D_h \dfrac{\partial^2\underline H}{\partial x^2}-
\dfrac{\beta_h \underline V (N^*_h-\underline H)}{N_h^*}+(d_h+\gamma_h)\underline H\\
& &\quad =\delta \psi[\lambda_0+\frac {\delta \beta_h\phi}{N^*_h}],
\end{eqnarray*}
for  $-h_0<x<h_0$ and $t>0$. Recalling $\lambda_0<0$, we can choose $\delta$ sufficiently small such that
\begin{eqnarray*}
\left\{
\begin{array}{lll}
\dfrac{\partial \underline V}{\partial t}-D_v \dfrac{\partial^2 \underline V}{\partial x^2}\leq
\dfrac{\beta_v(N_v^*-\underline V)\underline H}{N_h^*} -r_v(1-q) \underline V,\; & -h_0<x<h_0,\, t>0, \\
\dfrac{\partial \underline H}{\partial t}-D_h \dfrac{\partial^2\underline H}{\partial x^2}\leq
\dfrac{\beta_h \underline V (N^*_h-\underline H)}{N_h^*}-(d_h+\gamma_h)\underline H,\; &  -h_0<x<h_0, \, t>0,\\
\underline V(x,t)=\underline H(x, t)=0,&x=\pm h_0\, \, t>0,\\
0=-h'_0\geq -\mu D_h \frac{\partial \underline H}{\partial x}(-h_0, t), & t>0, \\
0=h'_0\leq -\mu D_h \frac{\partial \underline H}{\partial x}(h_0, t), & t>0,\\
\underline {V}(x,0)\leq V_{i0}(x),\ \underline{H}(x,0)\leq H_{i0}(x),\; &-h_0\leq x\leq h_0.
\end{array} \right.
\end{eqnarray*}
Hence by applying Lemma 3.2, we get that $V_i(x,t)\geq\underline V(x,t)$ and  $H_i(x,t)\geq\underline H(x,t)$
in $[-h_0, h_0]\times [0,\infty)$. It follows that $\displaystyle \liminf_{t\to
+\infty} \ |H_i(\cdot, t)||_{C([g(t), h(t)])}\geq \delta \psi(0)>0$ and therefore $h_\infty-g_\infty=+\infty$ by Lemma 5.2.

If $R_0^F(0):=R_0^D((-h_0, h_0))=1$, then for any positive time $t_0$, we have $g(t_0)<-h_0$ and $h(t_0)>h_0$, therefore $R_0^D((g(t_0), h(t_0)))>
R_0^D((-h_0, h_0))=1$ by the monotonicity in Theorem 3.1(iii).
Replacing the initial time $0$ by the positive time $t_0$, one can acquire $h_\infty-g_\infty=+\infty$ as above.
 \epf

\begin{rmk} It follows from the above proof that spreading happens if there exists $t_0\geq 0$ such that $R_0^F(t_0)\geq 1$.
\end{rmk}

Next, we consider the asymptotic behavior of the solution to problem \eqref{a3} when the spreading occurs.

\begin{thm} \label{asymp} If $R_0^F(t_0)(:=R_0^D((g(t_0), h(t_0))))\geq 1$ for some $t_0\geq 0$, then $h_\infty=-g_\infty=+\infty$ and
$$\lim_{t\to +\infty} \ (V_i(x,t),H_i(x,t))=(V_i^*, H_i^*)$$
uniformly in any bounded subset of $(-\infty, \infty)$, where $(V_i^*, H_i^*)$ is the unique positive equilibrium of
the corresponding ODE systems \eqref{a31}.
\end{thm}
\bpf We first present the limit superior of the solution.
It follows from the comparison principle that $(V_i(x,t),H_i(x,t))\leq (\overline V(t), \overline H(t))$
for $g(t)<x<h(t),\, t>0$, where
$(\overline V(t), \overline H(t))$ is the solution of the problem
\begin{eqnarray}
\label{ode1}
\left\{
\begin{array}{lll}
&\dfrac{\textrm{d} \overline V}{\textrm{d} t}=\dfrac{\beta_v(N_v^*-\overline V)\overline H}{N_h^*} -r_v(1-q)\overline V_i,\; &  t>0, \\
&\dfrac{\textrm{d} \overline H}{\textrm{d} t}=\dfrac{\beta_h \overline V (N^*_h-\overline H)}{N_h^*}-(d_h+\gamma_h)\overline H,\; & t>0,\\
&\overline V(0)=||V_{i0}||_{L^\infty([-h_0, h_0])},\ \overline H(0)=||H_{i0}||_{L^\infty([-h_0, h_0])}.&
\end{array} \right.
\end{eqnarray}
Since $R_0>R_0^F(t_0)\geq 1$, it is well known that the unique positive equilibrium $(V_i^*, H_i^*)$ is globally stable for the ODE system (\ref{ode1}) and $\displaystyle \lim_{t\to\infty}(\overline
V(t), \overline H(t))= (V_i^*, H_i^*)$. Therefore we deduce that
\begin{equation}\label{123}
\limsup_{t\to +\infty} \ (V_i(x,t), H_i(x,t))\leq (V_i^*, H_i^*)
\end{equation}
uniformly for $x\in (-\infty, \infty)$.

We now derive the limit inferior of the solution.
Thanks to $R_0>1$ by assumption, there is $L_0>0$ such that $R_0^D((-L_0, L_0))>1$.
Since $h_\infty=-g_\infty=+\infty$,
 for any $L\geq L_0$, there exists $t_L>0$ such that $g(t)\leq -L$ and $h(t)\geq L$ for $t\geq t_L$.
It follows from Lemma 4.2(b) that $R_0^F(t_{L_0})>1$, as in the proof of Theorem 6.1,
  we can choose $\delta$ sufficiently small such that $(V_i,H_i)\geq (\delta \phi, \delta \psi)$ in $[-L_0, L_0]\times [t_{L_0}, \infty)$,
which implies that the solution can not decay to zero.

We extend $\phi(x)$ to $\phi_{L_0}(x)$ by defining $\phi_{L_0}(x):=\phi(x)$ for $-L_0\leq x\leq L_0$ and $\phi_{L_0}(x):=0$ for $x<-L_0$ or $x>L_0$. Now for $L\geq L_0$, $(V_i, H_i)$ satisfies
\begin{eqnarray}
\left\{
\begin{array}{lll}
\dfrac{\partial V_i}{\partial t}-D_v \Delta V_i=\dfrac{\beta_v(N_v^*-V_i)H_i}{N_h^*} -r_v(1-q) V_i,\; & g(t)<x<h(t),\ t>t_L, \\
\dfrac{\partial H_i}{\partial t}-D_h \Delta H_i=\dfrac{\beta_h V_i (N^*_h-H_i)}{N_h^*}-(d_h+\gamma_h) H_i,\; &  g(t)<x<h(t),\ t>t_L,\\
V_i(x,t)=H_i(x,t)=0, \quad & x=g(t)\, \textrm{or}\, x=h(t),\ t>t_L,\\
V_i(x,t_L)\geq \delta \phi_{L_0},\ H_i(x,t_L)\geq \delta \psi_{L_0},
 & -L\leq x\leq L.
\end{array} \right.
\label{fs1}
\end{eqnarray}
Therefore, we have $(V_i, H_i)\geq (w, z)$ in $[-L, L]\times [t_L, \infty)$,
where $(w, z)$ satisfies
\begin{eqnarray}
\left\{
\begin{array}{lll}
w_t-D_v \Delta w=\dfrac{\beta_v(N_v^*-w)z}{N_h^*} -r_v(1-q) w,\; & -L<x<L,\ t>t_L, \\
z_t-D_h \Delta z=\dfrac{\beta_h (N^*_h-z)w}{N_h^*}-(d_h+\gamma_h) z,\; &  -L<x<L,\ t>t_L,\\
w(x,t)= z(x,t)=0, \quad & x=\pm L,\ t>t_L,\\
w(x,t_L)=\delta \phi_{L_0},\ z(x,t_L)=\delta \psi_{L_0},  & -L\leq x\leq L.
\end{array} \right.
\label{fs11}
\end{eqnarray}
The system (\ref{fs11}) is quasimonotone increasing. Therefore,
it follows from the upper and lower solution method
 and the theory of monotone dynamical systems ( \cite{HS}, Corollary 3.6) that
$\displaystyle \lim_{t\to +\infty} \ (w(x,t), v(x,t))\geq (w_L(x), z_L(x))$ uniformly in
$[-L, L]$, where $(w_L, z_L)$ satisfies
\begin{eqnarray} \label{fs12}\left\{
\begin{array}{lll}
-D_v \Delta w_L=\dfrac{\beta_v(N_v^*-w_L)z_L}{N_h^*} -r_v(1-q) w_L,\; &  -L<x<L,  \\
-D_h \Delta z_L=\dfrac{\beta_h (N^*_h-z_L)w_L}{N_h^*}-(d_h+\gamma_h) z_L, \; & -L<x<L,  \\
w_L(x)=z_L(x)=0, &x=\pm L
\end{array} \right.
\end{eqnarray}
and is the minimal solution over $(\delta \phi_{L_0},\ \delta \psi_{L_0})$.

The pair $(w_L(x), z_L(x))$ depends on $L$ and increases with $L$, that is, if $0<L_1<L_2$, then $(w_{L_1}(x), z_{L_1}(x))\leq (w_{L_2}(x),z_{L_2}(x))$
 in $[-L_1, L_1]$. The result is derived by comparing the boundary
conditions and initial conditions in (\ref{fs11}) for $L=L_1$ and $L=L_2$.

Let $L\to \infty$.  By classical elliptic regularity theory and a diagonal procedure,
it follows that  $(w_{L}(x), z_{L}(x))$ converges
uniformly on any compact subset of $(-\infty, \infty)$ to $(w_\infty, z_\infty)$,
which is continuous on $(-\infty, \infty)$ and satisfies
\begin{eqnarray*} \left\{
\begin{array}{lll}
-D_v \Delta w_{\infty}=\dfrac{\beta_v(N_v^*-w_{\infty})z_{\infty}}{N_h^*} -r_v(1-q) w_{\infty},\; &  -\infty<x<\infty,  \\
-D_h \Delta z_{\infty}=\dfrac{\beta_h (N^*_h-z_{\infty})w_{\infty}}{N_h^*}-(d_h+\gamma_h) z_{\infty} \; &  -\infty<x<\infty, \\
w_{\infty}(x)\geq \delta \phi_{L_0},\ z_{\infty}(x)\geq \delta \psi_{L_0}, &-\infty<x<\infty.
\end{array} \right.
\end{eqnarray*}

Next, we observe that $(w_\infty (x), z_{\infty}(x))\equiv (V_i^*, H_i^*)$, which can be derived by considering the corresponding reaction-diffusion system, whose solution tends to
the unique constant solution $(V_i^*, H_i^*)$.

Now for any given $[-M, M]$ with $M\geq L_0$, due to $(w_{L}(x), z_{L}(x))\to (V_i^*, H_i^*)$ uniformly in $[-M, M]$ as $L\to \infty$, we deduce that for any $\varepsilon >0$, there exists $L^*>L_0$ such that  $(w_{L^*}(x), z_{L^*}(x))\geq  (V_i^*-\varepsilon, H_i^*-\varepsilon)$ in $[-M, M]$. As above, there exists $t_{L^*}$ such that $[g(t), h(t)]\supseteq [-L^*, L^*]$ for $t\geq t_{L^*}$.
Therefore,
$$(V_i(x,t), H_i(x,t))\geq (w(x,t), z(x,t))\ \textrm{in}\ [-L^*, L^*]\times [t_{L^*}, \infty),$$ and
$$\lim_{t\to +\infty} \ (w(x,t), z(x,t))\geq (w_{L^*}(x), z_{L^*}(x))\ \textrm{in}\ [-L^*, L^*].$$
Using the fact that $(w_{L^*}(x), z_{L^*}(x))\geq (V_i^*-\varepsilon, H_i^*-\varepsilon)$ in $[-M, M]$ gives
 $$\liminf_{t\to +\infty} \ (V_i(x, t), H_i(x,t))\geq (V_i^*-\varepsilon, H_i^*-\varepsilon)\ \textrm{in}\ [-M, M].$$
Subsequently, the arbitrariness of $\varepsilon>0$ can lead to $\liminf_{t\to +\infty} \ H_i(x, t)\geq V_i^*$ and $\displaystyle \liminf_{t\to +\infty} \ H_i(x,t)\geq H_i^*$ uniformly in $[-M,M]$, which together with (\ref{123})
imply that $\displaystyle \lim_{t\to +\infty} \ V_i(x,t)=V_i^*$
and  $\displaystyle \lim_{t\to +\infty} \ H_i(x,t)=H_i^*$ uniformly in any bounded subset of $(-\infty, \infty)$.
\epf

\section{Basic reproduction numbers and spreading speeds}

In this paper, we have considered a simplified spatial model for West Nile virus (WNv)
which describes the diffusive transmission of the virus, and examined the dynamical behavior of the population $(V_i, H_i)$
with spreading fronts $x=g(t)$ and $x=h(t)$ determined by \eqref{a3}.  We have obtained several threshold conditions (basic reproduction numbers) and presented the asymptotic behavior results of the solutions. Sufficient conditions are given to ensure that the spreading or vanishing happens.
To conclude this paper, we present some remarks and discussion on the threshold basic reproduction numbers and estimates of spreading speed.

To characterize the temporal and spatial dynamics of the spreading, we defined
a threshold number $R^{F}_0(t)$
which was used to decide whether or not the spatial spread of the virus will happen, that is, if there exists a $t_0\geq 0$ such that
$R^{F}_0(t)\geq 1$,
then the virus is in a spreading scenario;
conversely, if $R_0<1$, the spreading or vanishing of the virus is decided
by the initial values and the expanding ability of the virus in the new area.

Now we have four basic reproduction numbers $R_0$, $R_0^N$, $R_0^D$ and $R_0^F(t)$ for the simplified spatial spreading model, which are
defined in homogenies environment, in a bounded environment with no flux,
in a bounded environment with hostile boundary and in a expanding environment, respectively. Each of the four cases
may happen due to the spatial variations, especially the geographical characteristics of the endemic area of the virus in Canada and USA.

If the environment is homogenous, the simplified WNv model can be described by the ODE system (\ref{a31}),
where the basic reproduction number is defined by
$$R_0=\ \sqrt{\frac {\beta_v\beta_hN_v^*}{r_v(1-q)N_h^*(d_h+\gamma_h)}}\ ,$$
see the definition and detailed calculations in \cite{VW} (see also page 10 of \cite{Lewis05}).

If the environment is bounded and heterogenous, the basic reproduction numbers $R_0^N$ and $R_0^D$ involve the spatial
 habitat characters.
 The properties of $R_0^N$ has been discussed in \cite{AL}
 for an SIS model, where the habitat $\Omega$ is characterized as {\bf high-risk}
( or {\bf low-risk} ) if the spatial average of the transmission rate is greater than ( or less than ) the spatial average of the recovery rate. It was shown that for low-risk domains, the disease-free equilibrium is stable ($R_0^D < 1$) if and only if the mobility of infected
individuals lies above a threshold value. For high-risk domains the disease-free
equilibrium is always unstable ($R_0^D> 1$). The properties of $R_0^N$ can been found in \cite{LLZ} for a Logistic model.

If the environment is heterogenous and expanding, we call the threshold value $R_0^F$ as {\bf the spatial-temporal risk index}, which involves not only habitat characteristic, but also the time with which the habitat is changing. It fully reflects the complexity of the spatial spreading of the virus, see also recent work in \cite{GKLZ, GLZ}.

Recently, to consider reaction-diffusion epidemic models with compartmental structure, the theory of the principal eigenvalue for an elliptic eigenvalue problem has been developed by Wang and Zhao in \cite{WZ2012}, the basic reproduction number was established by the spectral radius of a next infection operator and its computation formulae were given.
Furthermore, a reaction-diffusion model for Lyme disease with a spatially heterogeneous structure has been studied by Wang and Zhao in \cite{WZ2015},
 the basic reproduction number of the disease was established, see also \cite{YZ2016} for a nonlocal and time-delayed reaction-diffusion model for Lyme disease with a spatially heterogeneous structure.

For vector-borne diseases, another important issue for WNv is that once the virus starts spreading, no matter whether it will be in a scenario of spreading or vanishing,
 it is interesting and important to estimate the spreading speed.

To best understand the spreading speed, let us recall the following well-known results
for the diffusive logistic equation:
\begin{equation}
\label{tw}
 u_t-d\Delta u=u(a-bu),\; t>0,\; x\in \R.
\end{equation}

Define $c_{\min}$ by $\inf \{c\}$, there exists a traveling wave for $c$, which is called {\bf the minimal wave speed}.
Fisher \cite{Fi} and Kolmogorov et al
\cite{KPP} showed that $c_{\min} =2\sqrt{ad}$, that is, for any $|c|\geq c_{\min}$, there
exists a traveling wave solution $u(t,x):=W(x+ct)$ for \eqref{tw} with the property that
\[W'(y)>0 \mbox{ for } y\in \R^1,\; W(-\infty)=0,\;\; W(+\infty)=a/b;\]
no such solution exists if $|c|<c_{\min}$.

Define $c^*$ as {\bf the spreading speed} of a new population
$u(t,x)$ (governed by the above logistic equation) with which the region $\{x: \, u\backsim 1\}$ where the new species
dominates takes over the set $\{x: \, u(x,0)=0\}$ where the new species
is initially absent. It was shown in \cite{KPP} that $c^*=c_{\min} =2\sqrt{ad}$. Also see Section 4 in \cite{AW1}, it was proved that for such $u(t,x)$,
\[
\lim_{t\to\infty,\; |x|\leq (c^*-\epsilon)t}u(t,x)=0,\;\;\;\;
\lim_{t\to\infty,\;|x|\geq (c^*+\epsilon)t}u(t,x)=a/b\] for any small
$\epsilon>0$.

For the following diffusive logistic problem with free boundary,
\begin{eqnarray}
\left\{
\begin{array}{lll}
u_{t}-d u_{xx}=u(a-b u),\; & t>0, \ g(t)<x<h(t),  \\
u(t,g(t))=0,\; g'(t)=-\mu u_x(t, g(t)), \quad & t>0, \\
u(t, h(t))=0, \;h'(t)=-\mu u_x(t,h(t)), &t>0,\\
g(0)=-h_0,\; h(0)=h_0, \; u(0, x)=u_{0}(x),\; &-h_0\leq x\leq h_0,
\end{array} \right.
\label{f2}
\end{eqnarray}
when spreading happens, it was shown in \cite{DG, DL, GLL} that
 $\displaystyle \lim_{t\to +\infty} \frac {h(t)}t=\displaystyle \lim_{t\to +\infty} \frac {-g(t)}t=k_0$,
 where $k_0<c_{\min} =2\sqrt{ad}$ satisfies $\mu U'_{k_0}(0)=k_0$ and $U_{k_0}(y)$ satisfies
\begin{eqnarray}
\left\{
\begin{array}{lll}
-d U''+k_0 U'=a U-b U^2, \quad &y>0,  \\
\;\;\; U(0)=0. &
\end{array} \right.
\label{m9}
\end{eqnarray}
$k_0$ is regarded as {\bf the asymptotic spreading speed} of the free boundary problem.

 \begin{figure}[htbp] \centering
\includegraphics[width=12cm]{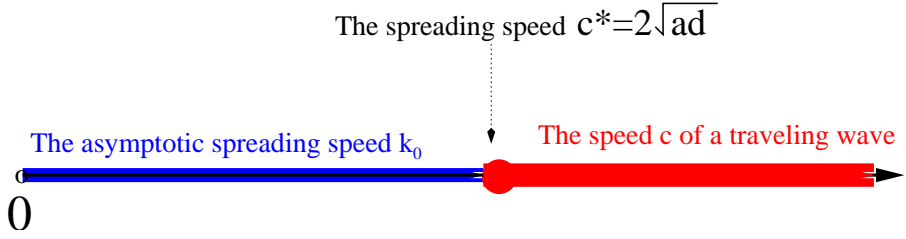}
   \caption{Comparison of the spreading speeds. $c^*$ is the spreading speed for the diffusive logistic equation $u_t-d\Delta u=u(a-bu)$, $c(\geq c_{\min}=2\sqrt{ad})$ is the speed of a traveling wave and $k_0(<2\sqrt{ad})$ is the asymptotic spreading speed of the corresponding free boundary problem. }\label{fig-speed}
\end{figure}

Notice that $U(0)=0$ and $\lim_{y\to \infty}U(y)=\frac ab$, it is actually the wavefront on half space and
therefore is called {\bf the semi-wavefront} with speed $k_0$ of the logistic equation,
which always exists uniquely for all $k_0<c_{\min}$.
It is well-known (\cite{DL}) that the asymptotic spreading speed of the free boundary problem (\ref{f2}) is a speed of the corresponding semi-wavefront, and the speed $k_0(\mu)$ is increasing with respect to $\mu$ and $k_0(\mu)\in (0, c_{\min})$ for $\mu\in (0, \infty)$. The relations of the spreading speeds are clearly shown in Fig. \ref{fig-speed}.

Now we look at the WNv models, the simplified spatial model reads
\begin{equation}\label{ab32}
\left \{
\begin{array}{lll}
\dfrac{\partial V_i}{\partial t}-D_v \dfrac{\partial^2 V_i}{\partial x^2}=\dfrac{\beta_v(N_v^*-V_i)H_i}{N_h^*} -r_v(1-q) V_i:=f_1(V_i,H_i),\; & x\in \R,\, t>0, \\
\dfrac{\partial H_i}{\partial t}-D_h \dfrac{\partial^2 H_i}{\partial x^2}=\dfrac{\beta_h V_i (N^*_h-H_i)}{N_h^*}-(d_h+\gamma_h) H_i:=f_2(V_i,H_i),\; & x\in \R, \, t>0.
\end{array} \right.
\end{equation}
As shown and discussed above, we have to ask the same questions:
\begin{itemize}
\item [1)] Does the minimal wave speed for the simplified spatial model \eqref{ab32} exist?

\item[2)]  Can the spreading speed be characterized as the minimal wave speed ($c^*=c_{\min}$)?

\item[3)]  Is the asymptotic spreading speed of the corresponding free boundary problem (\ref{a3})
less than the minimal wave speed ($k_0<c_{\min}$)?
\end{itemize}

The answer to the first question is certain. A traveling wave solution with speed $c$ for (\ref{ab32})
is a solution that possesses the form $(V_i(x+ct), H_i(x+ct))$
and the solution connects the disease-free and endemic equilibriums so that
$$\lim_{x+ct\to -\infty} (V_i, H_i)=(0,0)\ \textrm{and}\ \lim_{x+ct\to +\infty} (V_i, H_i)=(V^*_i, H^*_i).$$

We next adopt the theorem on existence of traveling waves provided in \cite{Li2005}.
It is easily verified that the reaction terms $\mathbf{f}$ (:=$(f_1,f_2)$) in (\ref{ab32}) satisfy the hypotheses of Theorem 4.2 in \cite{Li2005} with $\pmb{\beta}=(V^*_i, H^*_i)$. Using Theorem 4.2 in \cite{Li2005} (see also Theorem 4.1 in \cite{Lewis05}) yields the following traveling wave result.

\begin{thm} \label{trav} There exists a minimal speed $c_{\min}$ of traveling fronts such that for every $c\geq c_{\min}$
 the simplified diffusive system \eqref{ab32} admits a non-decreasing traveling wave solution
 $(V_i(x+ct), H_i(x+ct))$. If $c<c_{\min}$, no such traveling wave exists.
\end{thm}

As to the second question, the answer is also clear.
Since the only zeros of $\mathbf{f}$ are $\mathbf{0}$ and $(V^*_i, H^*_i)$, one can acquire from Theorem 4.2 in \cite{Li2005} (see also Theorem 4.2 in \cite{Lewis05}) the following conclusion.
\begin{thm} \label{trav1} The minimal speed of traveling fronts $c_{\min}$ for the simplified diffusive WNv model \eqref{ab32}
is equal to $c_0$, the spread speed for the system.
\end{thm}

For the calculation of the spread rate $c^*$, assuming the diffusion rate of the mosquitoes be small ($D_v\to 0$),
the spread rate for the simplified diffusive system \eqref{ab32} approaches the positive square root of the largest zero of a cubic, see Section 5 of \cite{Lewis05} in details.

To answer the third question 3), we need to consider the semi-wavefront of the simplified diffusive WNv model with free boundary (\ref{a3}),
\begin{equation}\label{ab33}
\left \{
\begin{array}{lll}
  -D_v \dfrac{\partial^2 V_i}{\partial x^2}+k_0 \dfrac{\partial V_i}{\partial x}=\dfrac{\beta_v(N_v^*-V_i)H_i}{N_h^*} -r_v(1-q) V_i,\; & x>0, \\
  -D_h \dfrac{\partial^2 H_i}{\partial x^2}+k_0 \dfrac{\partial H_i}{\partial x}=\dfrac{\beta_h V_i (N^*_h-H_i)}{N_h^*}-(d_h+\gamma_h) H_i,\; & x>0,\\
  V_i(0)=H_i(0)=0,\quad V_i(+\infty)=V_i^*,\, H_i(+\infty)=H_i^*,&
\end{array} \right.
\end{equation}
the existence of the solution, its properties and the answer to question 3) will be discussed in the future work.

To characterize the propagation of WNv, the existence of traveling waves was proved in \cite{Lewis05, Maidana09} and the spatial spread rate of infection in a reaction-diffusion model was calculated.
 A striking difference between the free boundary problem \eqref{a3} we have discussed here and the reaction-diffusion problem studied in \cite{Lewis05, Maidana09}
 is that the spreading fronts in \eqref{a3} are given explicitly by two functions
$x=g(t)$ and $x=h(t)$, the densities of infected birds and mosquitoes are $0$ beyond
the interval $(g(t), h(t))$; while in \cite{Lewis05, Maidana09},
the traveling wave solution is positive for all $x$
once $t$ is positive, which implies that infected birds and culex mosquitoes exist everywhere already. Moreover, the dynamics of \eqref{a3} exhibit a
spreading-vanishing dichotomy, which depends on the initial infected numbers and infected area. The reaction-diffusion model with free boundary
presents rich and complex dynamics about the spatial expanding of the infections
and we look forward to a further extension.

\end{document}